\documentclass[12pt]{amsart}

\usepackage{amssymb,verbatim,graphicx}
\usepackage[mathscr]{eucal}
\usepackage{enumerate}
\usepackage{xspace}
\usepackage{curves}

\newtheorem{theorem}{Theorem}[section]

\newtheorem{proposition}[theorem]{Proposition}
\newtheorem{corollary}[theorem]{Corollary}

\newtheorem{definition}[theorem]{Definition}

\theoremstyle{definition}
\newtheorem{example}[theorem]{Example}
\newtheorem{remark}[theorem]{Remark}

\newenvironment{cor-proof}[1]{\begin{trivlist} \item[] {\it Proof of Corollary #1.} }{\qed \end{trivlist} }
\newenvironment{thm-proof}[1]{\begin{trivlist} \item[] {\it Proof of Theorem #1.} }{\qed \end{trivlist} }

\numberwithin{equation}{section}
\numberwithin{figure}{section}

\def\&{\wedge}

\newcommand{\calL}{{\mathcal L}}

\newcommand{\E}{\mathbb{E}}
\newcommand{\R}{\mathbb{R}}

\newcommand{\bu}{\mathbf{u}}

\begin{document}

\title[Strings attached]{Strings attached: New light on an old problem
}

\author{Jeanne N. Clelland}
\address{Department of Mathematics, 395 UCB, University of
Colorado,
Boulder, CO 80309-0395}
\email{Jeanne.Clelland@colorado.edu}
\author{Peter J. Vassiliou}
\address{Program in Mathematics and Statistics, University of Canberra, Canberra, A.C.T., Australia, 2601}
\email{Peter.Vassiliou@canberra.edu.au}


\thanks{The first author was supported in part by NSF grant DMS-1206272.}

\begin{abstract}
The wave equation $u_{tt} = c^2 u_{xx}$ is generally regarded as a linear approximation to the equation describing the amplitude of a transversely vibrating elastic string in the plane.  But, as is shown in \cite{BC96}, the assumption of transverse vibration in fact implies that the wave equation describes the vibration precisely, with no need for approximation.  We give a simplified proof of this result, and we generalize to the case of an elastic string vibrating (transversely or not) in a Riemannian surface $M$.  In the more general setting, the assumption of transverse vibration is replaced by the assumption of ``perfect elasticity,'' and we show that the wave map equation $\nabla_{\bu_t} \bu_t = c^2 \nabla_{\bu_x} \bu_x$ gives a precise description of the vibration of a perfectly elastic string in $M$, with no need for approximation.  Finally, we give examples describing the motion of various vibrating strings in $\E^2$, $S^2$, and $\mathbb{H}^2$.
\end{abstract}

\maketitle

\section{Introduction}\label{intro-sec}

One of the most familiar partial differential equations in all of mathematics is the wave equation for a vibrating string,
\begin{equation}\label{1D-wave-eq}
u_{tt} = c^2 u_{xx},
\end{equation}
which describes the transverse vibration of an elastic string the in the plane.  The equilibrium position of the string is assumed to be an interval $[0, L]$ along the $x$-axis, and the function $u(x,t)$ represents the vertical displacement at time $t$ of the point on the string corresponding to the point with equilibrium position $(x,0)$.

A wide variety of methods have been used to derive equation \eqref{1D-wave-eq} from physical principles; we will critically review some of the best-known derivations (\cite{CH53}, \cite{CJ77}, \cite{Strauss92}) in \S \ref{history-sec}.  The methods discussed here vary in the specific techniques used during the derivation process, but they all make use of simplifying assumptions and/or approximations (e.g., constant tension, small amplitude vibrations), 
with the result that equation \eqref{1D-wave-eq} is deemed to be only a linear approximation to the ``true" wave equation.  Exceptions to this approach in the literature are rare; a few can be found in \cite{Antman80}, \cite{BC96}, \cite{Keller59}, \cite{Yong06}.

In \S \ref{better-derivation-sec}, we will start from scratch and carry out a derivation of the wave equation from physical principles, {\em without} any of the simplifying assumptions made in most derivations.  Our derivation is along the lines of that given in \cite{BC96}, but it is more straightforward.  Remarkably, it turns out that all the nonlinearities cancel each other out, and equation \eqref{1D-wave-eq} is the correct equation for the motion of a transversely vibrating elastic string---{\em not} ``just" a linear approximation.  In light of the simplicity of our derivation, we find it remarkable that this fact is not better known than it appears to be, and that the presumed ``nonlinear string equation" has persisted for many decades in textbooks.

In \S \ref{generalizations-sec}, we will consider generalizations of equation \eqref{1D-wave-eq} to the case of an elastic string with not-necessarily-transverse vibration, both in the Euclidean plane $\E^2$ and in an arbitrary Riemannian surface $M$, where the notion of ``transverse" vibration is not well-defined in general.  In this more general setting, the assumption of transverse vibration is replaced by the assumption of ``perfect elasticity" (cf. Definition \ref{perfect-elasticity-def}). We will see that with this assumption, the motion is governed by the wave map equation
\begin{equation}\label{wave-map-eq}
\nabla_{\bu_t} \bu_t = c^2 \nabla_{\bu_x} \bu_x
\end{equation}
for the function $\bu:[0,L] \times \R \to M$, where $\nabla$ represents the Levi-Civita connection of the Riemannian metric on $M$.  (We refer the reader to \cite{Tao12}, \cite{Struwe}, and \cite{TU04} for information on wave maps.)  Again, while the wave map is often portrayed as a ``first approximation" (see, e.g., \cite{Tao12}) to the motion of an elastic string, we will show that in fact no simplifying assumptions are needed in order to derive the wave map equation \eqref{wave-map-eq} from the physical principles governing the motion of the string.  In \S \ref{examples-sec} we give some examples of solutions to \eqref{wave-map-eq} in non-flat Riemannian surfaces.  Note that here the partial differential equation governing the motion of the string on a non-flat Riemannian surface is intrinsically nonlinear. This arises from the non-Euclidean nature of the ambient space in which the vibrations are taking place.

\section{Many derivations, many assumptions}\label{history-sec}

Consider an elastic string of constant linear density $\rho$, whose equilibrium position lies along the interval $[0, L]$ in the $x$-axis. Suppose that the string is allowed to vibrate only in the transverse direction, so that at time $t$, the position in the $xy$-plane of the point on the string corresponding to the point with equilibrium position $(x,0)$ is given by $(x,y) = (x, u(x,t))$.

\subsection{The Courant-Hilbert derivation}\label{CH-derivation-sec}

In Courant and Hilbert's classic text \cite{CH53}, the wave equation is derived from Lagrangian mechanics and the calculus of variations.  Simplifying assumptions include:
\begin{itemize}
\item The magnitude $T$ of the tension in the string is assumed to be constant.
\item The vibrations are assumed to have ``small" amplitude, so that the quantity $\sqrt{1 + u_x^2}$ can be approximated by its first-order Taylor polynomial $1 + \tfrac{1}{2} u_x^2$.
\end{itemize}
The kinetic energy $K$ of the string is given by
\[ K = \rho \int_0^L u_t^2\, dx, \]
while the potential energy $P$ is assumed to be proportional to the increase  in the string's length compared to its length at rest; i.e., 
\[ P = T \left( \int_0^L \sqrt{1 + u_x^2}\, dx - L \right) \approx \tfrac{1}{2} T \int_0^L u_x^2 \, dx. \]
This leads to the Lagrangian functional 
\[ \calL = K - P = \int_0^L \left(\rho u_t^2 - T u_x^2 \right)\, dx, \]
whose Euler-Lagrange equation is
\[ \rho u_{tt} - T u_{xx} = 0. \]
Setting $c^2 = \frac{T}{\rho}$ yields the wave equation \eqref{1D-wave-eq}.

\subsection{The Coulson-Jeffrey derivation}\label{CJ-derivation-sec}

Coulson and Jeffrey's well-known text \cite{CJ77} takes a different approach that avoids the calculus of variations, and instead considers the forces acting on a small segment of the string, corresponding to the interval $[x, x + \Delta x]$.
Simplifying assumptions include:
\begin{itemize}
\item The magnitude $T$ of the tension in the string is assumed to be constant.
\item While it is noted that the condition for $T$ to be approximately constant is that ``the wave disturbance is not too large," this assumption is not made explicit as in the Courant-Hilbert derivation until fairly late in the process.
\end{itemize}
The equation of motion for this segment of the string is given by Newton's law of motion $F=ma$; as the vibration is assumed to be transverse, only the vertical component of the motion is considered.  

\begin{figure}[h]
\setlength{\unitlength}{1.5pt}
\begin{center}
\begin{picture}(200,135)(0,20)
\thinlines
\put(20,30){\line(1,0){160}}
\qbezier(30,50)(80,130)(160,130)
\thicklines
\put(57,85.5){\vector(-1,-1){40}}
\thinlines
\put(57,85.5){\line(1,0){58}}
\put(57,85.5){\arc(10,0){45}}
\put(57,23){\makebox(0,0){$x$}}
\put(57,30){\line(0,1){55.5}}
\put(72,60){\makebox(0,0){$u(x,t)$}}
\put(85,94){\makebox(0,0){$\theta(x,t)$}}
\put(115,23){\makebox(0,0){$x + \Delta x$}}
\put(115,85.5){\line(0,1){37}}
\put(88,80){\makebox(0,0){$\Delta x$}}
\put(123,105){\makebox(0,0){$\Delta u$}}
\put(57,30){\circle*{2}}
\put(115,30){\circle*{2}}
\put(57,85.5){\circle*{2}}
\put(115,123.2){\circle*{2}}
\thicklines
\put(115,123.2){\vector(3,1){45}}
\thinlines
\put(30,75){\makebox(0,0){$T(x,t)$}}
\put(133,143){\makebox(0,0){$T(x+\Delta x,t)$}}
\put(115,123){\line(1,0){40}}
\put(115,123){\arc(10,0){20}}
\put(170,127){\makebox(0,0){$\theta(x+\Delta x,t)$}}
\end{picture} 
\end{center}
\caption{String vibrating transversely in the plane}
\label{classical-string-fig}
\end{figure}
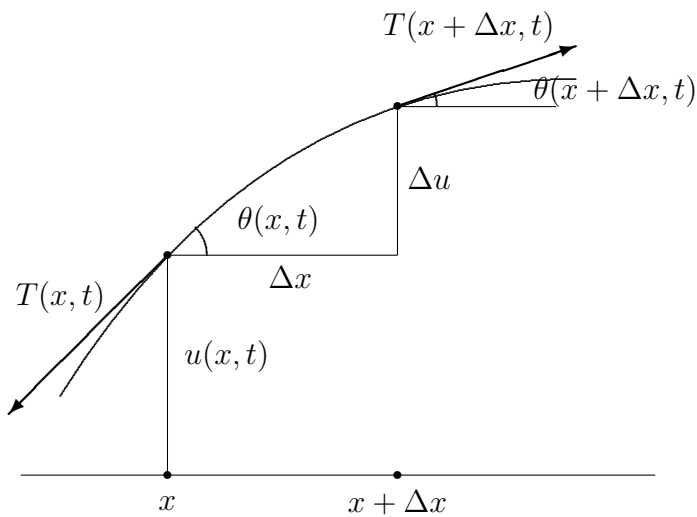

For each fixed $t$, let $\theta(x,t)$ denote the angle between the tangent vector to the string at the point $(x, u(x,t))$ and the horizontal, as shown in Figure \ref{classical-string-fig}. 

 The vertical component of the force exerted on the right-hand end of the string segment is given by $T \sin \theta(x + \Delta x,t)$, while the vertical component of the force exerted on the left-hand end of the string segment is given by $-T \sin \theta(x,t)$, so that the total vertical force acting on the string segment is
\[ F = T \sin \theta(x + \Delta x, t) - T \sin \theta(x, t) . \]
The arc length of the string segment is approximately 
\[ \Delta s = \sqrt{1 + u_x^2}\, \Delta x, \]
and so the mass of the string segment is assumed to be 
\[ m = \rho \Delta s = \rho \sqrt{1 + u_x^2}\, \Delta x. \]
Note that this assumption is incorrect: the mass of this string segment at rest is equal to $\rho \Delta x$, not $\rho \Delta s$, and conservation of mass implies that the mass remains equal to $\rho \Delta x$ regardless of how the string is stretched.  (This error appears frequently---although not universally---in the literature; it is generally inadvertently remedied at a later stage in the derivation where $u_x^2$ is assumed to be small, and therefore $\Delta s \approx \Delta x$.) Since $a = u_{tt}$, Newton's equation becomes:
\[ \rho \, u_{tt} \sqrt{1 + u_x^2}\,  \Delta x = T \sin \theta(x + \Delta x, t) - T \sin \theta(x, t). \]
Dividing by $\Delta x$ and taking the limit as $\Delta x \to 0$ yields
\begin{align}
 \rho \, u_{tt} \sqrt{1 + u_x^2} & = T \frac{\partial}{\partial x} (\sin \theta(x, t)) \label{CJ-derivation-eq1} \\
 & = T (\cos \theta(x, t)) \theta_x. \notag
\end{align}
The trigonometric identities
\begin{equation}
 \tan \theta = u_x, \qquad \sin \theta = \frac{u_x}{\sqrt{1 + u_x^2}}, \qquad  \cos \theta = \frac{1}{\sqrt{1 + u_x^2}} \label{trig-identities-1}
\end{equation}
imply that
\[ u_{xx} = \frac{\partial}{\partial x} (\tan \theta) = (\sec^2 \theta)\theta_x = (1 + u_x^2) \theta_x. \]
Therefore,
\begin{equation}
 \theta_x = \frac{u_{xx}}{(1 + u_x^2)}, \label{trig-identities-2}
\end{equation}
and equation \eqref{CJ-derivation-eq1} is equivalent to the {\em nonlinear} equation
\begin{equation}
 \rho \, u_{tt} = \frac{T u_{xx}}{(1 + u_x^2)^2}. \label{CJ-derivation-eq2}
\end{equation}
(If the mass of the string segment had been correctly represented as $\rho \Delta x$ instead of $\rho \Delta s$, equation \eqref{CJ-derivation-eq2} would instead have become
\begin{equation}
 \rho \, u_{tt} = \frac{T u_{xx}}{(1 + u_x^2)^{3/2}}.) \label{CJ-derivation-eq3}
\end{equation}
Only at this stage in the derivation do Coulson and Jeffrey apply the assumption that $u_x^2$ is small, so that equation \eqref{CJ-derivation-eq2} may be approximated by the linearized form \eqref{1D-wave-eq}, with $c^2 = \frac{T}{\rho}$.

\subsection{The undergraduate version}\label{ugrad-derivation-sec}

Yet another variant on this derivation appears in many undergraduate PDE textbooks; as an example, we present a slight variation on the derivation given in \cite{Strauss92}.  In this version, assumptions about the tension $T$ and the magnitude of $u_x$ are postponed until later in the process.  As in the Coulson-Jeffrey derivation, this argument is based on Newton's equation applied to the string segment corresponding to the interval $[x, x + \Delta x]$, but this time Newton's equation appears in its vector form $\mathbf{F} = m \mathbf{a}$.  With the same notation as in \S \ref{CJ-derivation-sec} (and using $m = \rho \Delta x)$, the vertical component of Newton's equation becomes
\begin{equation}
\rho \, u_{tt} \Delta x = T(x + \Delta x, t) \sin \theta(x + \Delta x, t) - T(x, t) \sin \theta(x, t), \label{ugrad-derivation-eq-vert}
\end{equation}
while the horizontal component becomes (due to the assumption of transversality)
\begin{equation}
 T(x + \Delta x, t) \cos \theta(x + \Delta x, t) - T(x, t) \cos \theta(x, t) = 0. \label{ugrad-derivation-eq-horiz}
\end{equation}
First consider equation \eqref{ugrad-derivation-eq-horiz}.  Dividing by $\Delta x$ and taking the limit as $\Delta x \to 0$ yields
\[ \frac{\partial}{\partial x} \left(T(x,t) \cos \theta(x,t) \right) =
T_x (\cos \theta) - T (\sin \theta) \theta_x = 0. \]
{\em Now} the simplifying assumption is made that $\theta$ is ``small," and therefore 
\[ \sin \theta \approx 0, \qquad \cos \theta \approx 1; \]
hence, $T_x \approx 0$.  So $T$ is now assumed to be approximately constant, but in this case the assumption of constant tension is a {\em consequence} of the hypothesis of transverse motion and the assumption that the vibrations are ``small," rather than a separate assumption. 

With these assumptions in hand, equation \eqref{ugrad-derivation-eq-vert} becomes
\[ \rho \, u_{tt} \Delta x  = T\sin \theta(x + \Delta x, t) - T \sin \theta(x, t), \]
and essentially the same argument as in \S \ref{CJ-derivation-sec} shows that, with the assumption that $u_x^2 \approx 0$, this equation is approximately equivalent to equation \eqref{1D-wave-eq}.

\section{A more accurate derivation}\label{better-derivation-sec}

In this section we present an alternative derivation for equation \eqref{1D-wave-eq}.  We will follow the general strategy of the derivation of \S \ref{ugrad-derivation-sec}, but we will assume {\em only} that the vibration is transverse, so that motion occurs only in the vertical direction.
No a priori assumptions will be made about either the function $T(x,t)$ (except that it is assumed to be differentiable) or the magnitude of the vibrations.  (Derivations using only these assumptions are given in \cite{BC96} and \cite{Keller59}; ours is similar, but more straightforward.)

As in \S \ref{ugrad-derivation-sec}, we apply Newton's equation $\mathbf{F} = m \mathbf{a}$ to the string segment corresponding to the interval $[x, x + \Delta x]$, which leads to equations \eqref{ugrad-derivation-eq-vert} and \eqref{ugrad-derivation-eq-horiz}.  Dividing both of these equations by $\Delta x$ and taking the limit as $\Delta x \to 0$ yields the equations
\begin{gather}
\rho\, u_{tt} = \frac{\partial}{\partial x} \left( T(x,t) \sin \theta(x,t) \right),   \label{right-derivation-eq-vert} \\
\frac{\partial}{\partial x} \left( T(x,t) \cos \theta(x,t) \right) = 0. \label{right-derivation-eq-horiz}
\end{gather}
First consider equation \eqref{right-derivation-eq-horiz}.  Write this equation as
\[ T_x (\cos \theta) - T (\sin \theta) \theta_x = 0, \]
and substitute in the trigonometric identities \eqref{trig-identities-1}, \eqref{trig-identities-2} to obtain
\begin{equation}
 T_x \frac{1}{\sqrt{1 + u_x^2}} - T \frac{u_x u_{xx}}{(1 + u_x^2)^{3/2}} = 0.   \label{right-derivation-eq-3}
\end{equation}
Equation \eqref{right-derivation-eq-3} may be regarded as a separable differential differential equation for the function $T(x,t)$: rewrite equation \eqref{right-derivation-eq-3} as
\[ \frac{T_x}{T} = \frac{u_x u_{xx}}{(1 + u_x^2)}. \]
Now integrate and exponentiate to obtain:
\begin{equation}\label{right-derivation-T-eq}
T(x,t) = C(t) \sqrt{1 + u_x^2},
\end{equation}
where $C(t)$ is an arbitrary function of $t$.  

Observe that, for each fixed $t$, equation \eqref{right-derivation-T-eq} says that the magnitude of the tension of the string at the point $(x, u(x,t))$ is proportional to the derivative of the arc length function 
\[ s(x,t) = \int_0^x \sqrt{1 + u_x(\chi, t)}\, d\chi \]
with respect to the curve parameter $x$.
Assuming that the physical parameters of the string are not changing in time, it is reasonable to assume that the ``constant" of proportionality $C(t)$ is, in fact, equal to a constant $T_0$, which represents the tension of the string in its equilibrium position.  The quantity 
\[ \sigma = \frac{\partial s}{\partial x} = \sqrt{1 + u_x^2} \]
may be regarded as the ``stretching factor" of the string at the point $(x, u(x,t))$, and, following \cite{Yong06}, we make the following definition:

\begin{definition}\label{perfect-elasticity-def}
A string is called {\em perfectly elastic} provided that, when the string at tension $T_0$ is stretched by a factor of $\sigma$, the tension becomes equal to $T = T_0 \sigma$.
\end{definition}

Thus, equation \eqref{right-derivation-T-eq} may be interpreted as saying that the hypothesis of transverse vibrations {\em implies} that the string is perfectly elastic, merely as a consequence of Newton's law.  This observation is significant enough that we state it as:

\begin{proposition}\label{perfect-elasticity-prop}
The tension of an elastic string vibrating transversely in the plane must satisfy the perfect elasticity condition \eqref{right-derivation-T-eq}.
\end{proposition}

\begin{remark}
Perfect elasticity is a special---and idealized---case of the more familiar notion of ``linear elasticity."  Linear elasticity refers to a string that satisfies Hooke's Law, which in its infinitesimal version says that, at a point where the string's stretch factor is equal to $\sigma$, the tension in the string is proportional to $(\sigma - \sigma_0)$, where $\sigma_0$ is the value of $\frac{\partial s}{\partial x}$ for which the tension is equal to zero.  Perfect elasticity is the special case where $\sigma_0=0$, which mathematically corresponds to a string which shrinks to a point in the absence of external forces.  As a practical matter, this condition is almost certainly never precisely satisfied by real physical strings; however, a rubber band or a spring stretched to several times its natural length may satisfy it to a reasonable approximation.  (Physical data supporting this assertion may be found in \cite{Keller59}.)  Consequently, the vibrating motion of real physical springs can never be truly transverse, but only approximately so.
\end{remark}

Continuing with the derivation, substitute the trigonometric identity \eqref{trig-identities-1} and the expression \eqref{right-derivation-T-eq} for $T(x,t)$ into equation \eqref{right-derivation-eq-vert} to obtain
\begin{align}
 \rho \, u_{tt} & = \frac{\partial}{\partial x} \left( C(t) \sqrt{1 + u_x^2} \frac{u_x}{\sqrt{1 + u_x^2}} \right ) \notag \\
 & =  \frac{\partial}{\partial x} \left( C(t) u_x \right) \label{right-derivation-final-answer} \\[0.1in]
 & = C(t) u_{xx}. \notag
\end{align}
With the assumption that $C(t) = T_0$, equation \eqref{right-derivation-final-answer} is equivalent to \eqref{1D-wave-eq}, with $c^2 = \frac{T_0}{\rho}$.

\section{Vibrating strings in Riemannian surfaces}\label{generalizations-sec}

Now consider the problem of an elastic string of constant linear density $\rho$ vibrating in a Riemannian surface $M$.  Even the innocuous-sounding assumption of ``constant linear density" merits some examination.  If a string is represented by a parametrized curve
\[ \bu: I \to M, \]
for some interval $I \subset \R$, then the {\em density} of the string is represented by a 1-form $d\mu$ on $I$, defined by the condition that the mass of any segment $\bu([a,b])$ of the string is given by
\begin{equation*}
 m(a,b) = \int_a^b d\mu. \label{define-mass-eqn}
\end{equation*}
``Constant linear density" is generally taken to mean that $d\mu = \rho\, dx$ for some constant $\rho$, but in fact this condition is dependent on the parametrization of the string.  A string with any smooth, nonvanishing density $1$-form $d\mu$ can always be given a reparametrizion $\bu(m)$ with the property that $d\mu = \rho\, dm$ for some constant $\rho>0$.  
(In fact, by taking $m$ to be the function 
\[ m(x) = \int_0^x d\mu, \]
we can arrange that $\rho=1$; however, it is often convenient to allow other values of $\rho$, so we will allow $\rho$ to be any positive constant.)  
  
\begin{definition}
Let $I \subset \R$ be a connected interval, and let $m$ be a local coordinate on $I$.  A parametrized curve $\bu:I \to M$ in a Riemannian surface $M$ with a specified density 1-form $d\mu$ on $I$ is said to be {\em parametrized by constant density} if $d\mu = \rho\, dm$ for some constant $\rho >0$.
\end{definition}

As observed above, a string with any smooth, nonvanishing density 1-form can be pa\-ra\-me\-trized by density, so the assumption of ``constant linear density" must mean something more than simply having the string parametrized by constant density.  In the classical case of a string vibrating transversely in the plane, the condition of constant linear density is characterized by the fact that the {\em equilibrium} position of the string has a parametrization which is simultaneously a parametrization by constant density {\em and} a parametrization by arc length.  As the vibrating string expands and contracts, however, the parametrization $(x, u(x,t))$ for the position of the string at a given time $t$ generally does not remain an arc length parametrization for the curve---but it {\em does} remain a constant density parametrization.

Next, consider the notion of ``transverse vibrations."  In the classical case of the string vibrating in the plane, ``transverse" is taken to mean ``in the vertical direction."  This notion relies crucially on the canonical local coordinates $(x,y)$ on the plane, and on the vertical vector field $\frac{\partial}{\partial y}$ defined on any open neighborhood in the plane.  But for a general Riemannian surface, there are no canonical local coordinates, and hence there is no obvious way to define the notion of a ``transverse" vibration.  (Even if we attempted to accomplish this by defining a ``transverse" vector field orthogonal to a curve representing the equilibrium position of some particular string, there is still no canonical way to extend such a vector field to an open neighborhood of the given curve segment.)  

Therefore, in light of Proposition \ref{perfect-elasticity-prop} we replace the assumption of ``transverse vibration" with an assumption of perfect elasticity, which does make sense even in the more general setting of a Riemannian surface.  But this too requires some care: according to Definition \ref{perfect-elasticity-def}, a string is perfectly elastic if the tension in the string is directly proportional to the stretch factor $\frac{\partial s}{\partial x}$.  But like the constant linear density condition, this condition is dependent on the parametrization of the string.  Since Definition \ref{perfect-elasticity-def} was formulated under the assumption of a constant density parametrization for the string, we now make the following refinement in order to make this assumption explicit:

\begin{definition}\label{perfect-elasticity-def-2}
A string in a Riemannian surface $M$ with a constant density  parametrization $\bu(m,t)$ is called {\em perfectly elastic} for $\sigma = \frac{\partial s}{\partial m}$ in the range $a \leq \sigma \leq b$ provided that, when the string at tension $T_0$ is stretched by a factor of $\sigma \in [a,b]$, the tension becomes equal to $T = T_0 \sigma$.
\end{definition}

In order to derive the equation of motion for the general case of a string vibrating in a Riemannian surface $M$, we will assume that:
\begin{itemize}
\item the string is parametrized by constant density;
\item the string is perfectly elastic;
\item the motion of the string is governed by Newton's Law.
\end{itemize}

Corollary \ref{constant-density-cor} will show that these conditions suffice to guarantee that a constant density parametrization for a string in equilibrium position is also a constant speed parametrization with respect to an arc length parameter along the curve segment.  This is the natural analog to the ``constant linear density" assumption in the classical case; the fact that it follows from the assumption of perfect elasticity suggests that perhaps the notion of perfect elasticity is more subtle than it might appear.

We are now ready to state our main result:

\begin{theorem}\label{main-theorem}
Let $M$ be a Riemannian surface with Levi-Civita connection $\nabla$, and let $\bu:[0,L] \times \R \to M$ be a smooth map, where for each fixed $t \in R$, the curve $\bu([0,L] \times \{t\})$ represents the position of a string vibrating in $M$ at time $t$.  Suppose that the string is parametrized by constant density and perfectly elastic, and that its motion is governed by Newton's Law.  Then the map $\bu(m,t)$ satisfies the wave map equation \begin{equation}\label{wave-map-eq2}
\nabla_{\bu_t} \bu_t = c^2 \nabla_{\bu_m} \bu_m.
\end{equation}
\end{theorem}

\begin{corollary}\label{constant-density-cor}
The equilibrium position of a string in a Riemannian surface $M$ satisfying the conditions of Theorem \ref{main-theorem} consists of a segment of a geodesic curve in $M$.  Furthermore, the given constant density parametrization of this curve segment is also a constant-speed parametrization with respect to an arc length parameter along the curve segment.
\end{corollary}

\begin{cor-proof}{\ref{constant-density-cor}}
Any equilibrium solution is independent of $t$; therefore $\nabla_{\bu_t} \bu_t = 0$.  Then \eqref{wave-map-eq2} implies that the solution curve $\bu(m,t) = \bu(m)$ satisfies the geodesic equation
\[ \nabla_{\bu_m} \bu_m = 0, \]
which implies that $\bu(m)$ is a constant-speed geodesic segment in $M$.
\end{cor-proof}

For the sake of clarity, we will prove Theorem \ref{main-theorem} in two steps: first we give a proof in the case that $M$ is the flat plane $\E^2$, where we can take advantage of the canonical identification between tangent planes at each point of $\E^2$.  Then we will indicate how the proof may be generalized for an arbitrary Riemannian surface $M$.

\begin{thm-proof}{\ref{main-theorem}}
First suppose that $M = \E^2$, and let 
\[ \bu(m,t) = \left( x(m,t), y(m,t) \right). \]
For fixed $t$, consider the string segment corresponding to the interval $[m, m + \Delta m]$ and the (vector-valued) tension forces $\mathbf{T}(m,t), \mathbf{T}(m + \Delta m, t)$ acting at the endpoints of the string segment.  (See Figure \ref{plane-string-fig}.)
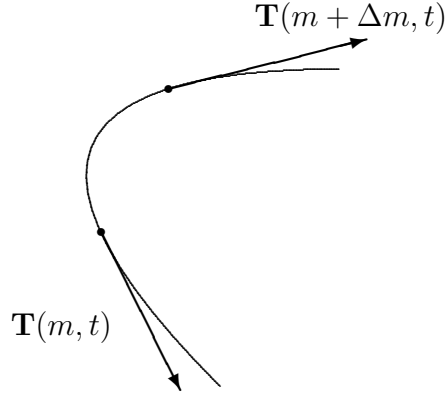
\begin{figure}[h]
\setlength{\unitlength}{1.5pt}
\begin{center}
\begin{picture}(200,135)(0,20)
\thinlines
\qbezier(100,50)(20,130)(130,130)
\thicklines
\put(70,89){\vector(1,-2){20}}
\thinlines
\put(70,89){\circle*{2}}
\put(87,125){\circle*{2}}
\thicklines
\put(87,125){\vector(4,1){50}}
\thinlines
\put(60,65){\makebox(0,0){$\mathbf{T}(m,t)$}}
\put(133,143){\makebox(0,0){$\mathbf{T}(m+\Delta m,t)$}}
\end{picture} 
\end{center}
\caption{Perfectly elastic string vibrating in the plane}
\label{plane-string-fig}
\end{figure}
According to the assumption of perfect elasticity, the magnitude of the force vector $\mathbf{T}(m,t)$ is equal to 
\[ T(m,t) = T_0\, \frac{\partial s}{\partial m}(m,t) = T_0 \left( \sqrt{x_m^2 + y_m^2} \right) \Big{\vert}_{(m,t)}\]
for some positive constant $T_0$; similarly, the magnitude of the force vector $\mathbf{T}(m+\Delta m,t)$ is equal to 
\[ T(m+\Delta m,t) =  T_0\, \frac{\partial s}{\partial m}(m + \Delta m,t) = T_0 \left( \sqrt{x_m^2 + y_m^2} \right) \Big{\vert}_{(m+\Delta m,t)}. \]
The direction of $\mathbf{T}(m + \Delta m,t)$ is given by the unit tangent vector 
\[ \mathbf{t}(m+\Delta m,t) = \left( \frac{1}{\sqrt{x_m^2 + y_m^2}} \left( x_m, y_m \right) \right) \Bigg{\vert}_{(m+\Delta m,t)} \] 
to the curve at the point $\mathbf{u}(m+\Delta m,t)$, while the direction of $\mathbf{T}(m,t)$ at the point $\mathbf{u}(m,t)$ is given by  
\[ -\mathbf{t}(m,t) = \left( \frac{-1}{\sqrt{x_m^2 + y_m^2}} \left( x_m, y_m \right) \right) \Bigg{\vert}_{(m,t)}. \] 
Therefore, we have
\begin{align*}
 \mathbf{T}(m,t) & = -T(m,t) \mathbf{t}(m,t) = -\left( T_0\,\sqrt{x_m^2 + y_m^2}\right)  \left( \frac{1}{\sqrt{x_m^2 + y_m^2}}
  (x_m, y_m) \right) \Bigg{\vert}_{(m,t)} \\
  &  = -T_0\, \mathbf{u}_m(m,t), \\[0.1in]
 \mathbf{T}(m+\Delta m,t) & = T(m+\Delta m,t) \mathbf{t}(m+\Delta m,t) = \left( T_0\,\sqrt{x_m^2 + y_m^2}\right)  \left( \frac{1}{\sqrt{x_m^2 + y_m^2}}
  (x_m, y_m) \right) \Bigg{\vert}_{(m+\Delta m,t)} \\
  &  = T_0\, \mathbf{u}_m(m+\Delta m,t) .
\end{align*}
Thus the total force acting on the string segment is
\begin{equation}
 \mathbf{F} = \mathbf{T}(m,t) + \mathbf{T}(m+\Delta m,t) = T_0 \left(\mathbf{u}_m(m+\Delta m,t)  -  \mathbf{u}_m(m,t)\right). \label{total-force-eqn}
\end{equation}
Since $\bu(m,t)$ is assumed to be a constant density parametrization, the mass of the string segment is given by $\rho\, \Delta m$.  Now Newton's equation 
becomes
\begin{equation}
 (\rho\, \Delta m) \, \mathbf{u}_{tt} = T_0 \left[\mathbf{u}_m(m+\Delta m,t)  -  \mathbf{u}_m(m,t)\right]. \label{take-my-limit}
\end{equation}
Dividing by $\Delta m$ and taking the limit as $\Delta m \to 0$ yields 
\[ \rho\, \mathbf{u}_{tt} = T_0\, \mathbf{u}_{mm}. \]
Taking $c^2 = \frac{T_0}{\rho}$ yields
\[  \mathbf{u}_{tt} = c^2\,  \mathbf{u}_{mm}, \]
which is equivalent to \eqref{wave-map-eq2} for the standard flat connection $\nabla$ on $\E^2$.

For a general Riemannian surface $M$, we must make the following modifications to this proof:
\begin{itemize}
\item In Newton's equation, the acceleration vector $\bu_{tt}$ must be replaced by $\nabla_{\bu_t} \bu_t$.
\item The expression
\[    \mathbf{T}(m,t) + \mathbf{T}(m+\Delta m,t) \]
in \eqref{total-force-eqn} no longer makes sense because the tangent vectors $\mathbf{T}(m,t)$, $\mathbf{T}(m+\Delta m,t)$ are based at different points of $M$.  In order to remedy this, we must replace the vector $\mathbf{T}(m+\Delta m,t)$ by its parallel transport backwards along the curve $\bu(\cdot,t)$ from $\bu(m+\Delta m, t)$ to $\bu(m,t)$.  
\end{itemize}
Substituting these changes into equation \eqref{take-my-limit}, dividing by $\Delta m$ and taking the limit as $\Delta m \to 0$ yields the wave map equation \eqref{wave-map-eq2}.

\end{thm-proof}

\section{Examples}\label{examples-sec}

In this section we compute several examples of solutions to equation \eqref{wave-map-eq2} (with $c=1$) in various Riemannian surfaces.  These examples illustrate how solutions may be affected not only by the geometry of the surface $M$, but also by the choice of a density 1-form $d\mu$ along the string.  (The choice of $d\mu$ is implicit in the choice of parametrization for the initial data curve, as this curve is assumed to be parametrized by constant density.)  

For each example we will specify a Riemannian surface $M$ and an initial curve segment $\bu(m, 0) = \bu_0(m)$ in $M$.  We will give an explicit representation for the PDE system \eqref{wave-map-eq2} in terms of local coordinates on $M$, and we will numerically solve this system for the given initial data, assuming fixed endpoints and zero initial velocity.\footnote{All numerical solutions were computed in Maple 16 using the ``pdsolve/numeric" function with standard options.}

\begin{example}\label{classic-wave-in-plane-ex} 
Let $M = \E^2$ be the flat plane.  With local coordinates $(x,y)$ and the standard flat metric on $\E^2$, the wave map equation \eqref{wave-map-eq2} is equivalent to the system
\begin{equation}
x_{tt} = x_{mm}, \qquad y_{tt} = y_{mm} \label{R2-system}
\end{equation}
for the functions $x(m,t), y(m,t)$.  

Suppose that the initial curve is the graph of $y = \sin (\pi x)$ for $0 \leq x \leq 1$.   The classical case of ``constant linear density" along this curve is represented by the parametrization
\[ x(m,0) = m, \qquad y(m,0) = \sin(\pi m), \qquad 0 \leq m \leq 1. \]
This parametrization leads to the classic solution
\[ x(m,t) = m, \qquad y(m,t) = \cos(\pi t) \sin(\pi m). \]
Some curves in this evolution are shown in Figure \ref{classic-string-in-R2-fig}; the initial curve is drawn as a thick curve, while subsequent curves in the evolution are thinner.
\begin{figure}[h]
\includegraphics[width=2.5in]{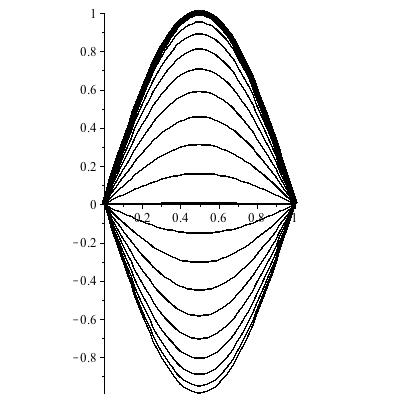}
\caption{Evolution of ``constant linear density" string in $\E^2$}
\label{classic-string-in-R2-fig}
\end{figure}

Now suppose that instead of the density 1-form $d\mu = dx$, we choose $d\mu = 2(x+\tfrac{1}{2})\, dx$, so that the density of the string is an affine linear function of the initial $x$-coordinate of the string.  This corresponds to choosing $m = (x + \tfrac{1}{2})^2$, which yields the parametrization 
\[ x(m, 0) = \sqrt{m} - \tfrac{1}{2}, \qquad y(m,0) = \sin(\pi (\sqrt{m} - \tfrac{1}{2})), \qquad \tfrac{1}{4} \leq m \leq \tfrac{9}{4} \]
for the initial data curve.
The system \eqref{R2-system} can be solved explicitly, either by Fourier series or by D'Alembert's formula; some curves in this evolution are shown in Figure \ref{modified-string-in-R2-fig}.
\begin{figure}[h]
\includegraphics[width=2.5in]{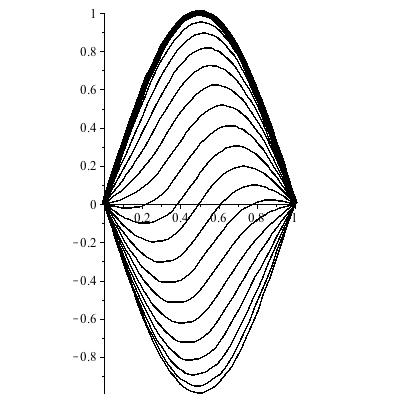}
\caption{Evolution of string with density $d\mu = 2(x+\tfrac{1}{2})\, dx$ in $\E^2$}
\label{modified-string-in-R2-fig}
\end{figure}
Note that, regardless of the initial curve, the general theory of the linear wave equation guarantees that the fixed-endpoint problem will have periodic solutions.  The period, however, depends on the parametrization: if the initial curve has parametrization
\[ x(m,0) = x_0(m), \qquad y(m,0) = y_0(m), \qquad a \leq m \leq b, \]
then the vibration of the string will have period equal to $2(b-a)$.

\end{example}

\begin{example}

Let $M = S^2$ be the unit sphere in $\R^3$, with local parametrization
\[ (x,y) \to (\cos x \cos y, \sin x \cos y, \sin y), \qquad 0 \leq x \leq 2\pi, \ \ \ -\frac{\pi}{2} < y < \frac{\pi}{2}. \]
The standard metric on $S^2$ is given by
\[ ds^2 = (\cos^2 y) \, dx^2 + dy^2, \]
and the wave map equation on $S^2$ is equivalent to the system
\begin{align}
x_{tt} & = x_{mm} - 2\tan y\, (x_m y_m - x_t y_t), \label{S2-system} \\
y_{tt} & = y_{mm} + \sin y \cos y \,(x_m^2 - x_t^2). \notag
\end{align}
This is a nonlinear system, and its solutions are not periodic in general.

By analogy with our previous example, we will consider ``sinusoidal" initial curves.  First suppose that the initial curve is described in local coordinates on $S^2$ by the equation
\[ y =  \frac{1}{2} \sin \left(\frac{2\pi}{3} x \right), \qquad 0 \leq x \leq \frac{3}{2}, \]
with parametrization
\[ x(m,0) = m, \qquad y(m,0) = \frac{1}{2} \sin \left(\frac{2\pi}{3} m \right), \qquad 0 \leq m \leq \frac{3}{2}. \]
Figure \ref{string-in-sphere-fig-1} shows the evolution of the curve through three cycles in the vibration; as in the last example, the thick curve in each picture represents the initial curve for that portion of the evolution.  Note that this solution is only quasi-periodic, as the shape of the wave changes with each successive cycle.
\begin{figure}[h]
\includegraphics[width=1.8in]{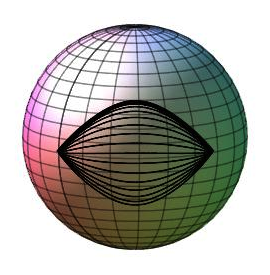}
\includegraphics[width=1.8in]{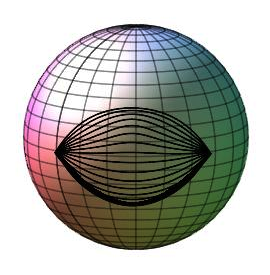}
\includegraphics[width=1.8in]{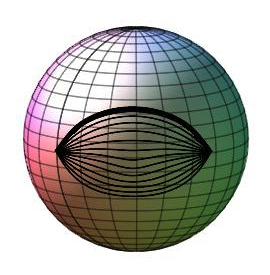}
\includegraphics[width=1.8in]{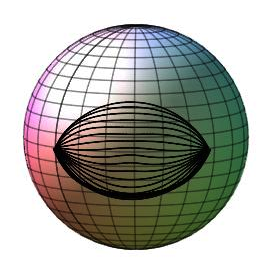}
\includegraphics[width=1.8in]{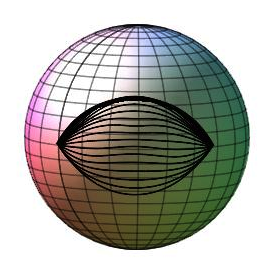}
\includegraphics[width=1.8in]{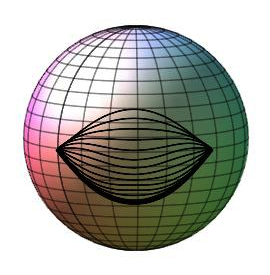}
\caption{Evolution of small amplitude string in $S^2$}
\label{string-in-sphere-fig-1}
\end{figure}

The effect of the nonlinearity becomes more significant if we increase the amplitude of the initial curve: suppose that the initial curve is described in local coordinates on $S^2$ by the equation
\[ y =  \sin \left(\frac{2\pi}{3} x \right), \qquad 0 \leq x \leq \frac{3}{2}, \]
with parametrization
\[ x(m,0) = m, \qquad y(m,0) = \sin \left(\frac{2\pi}{3} m \right), \qquad 0 \leq m \leq \frac{3}{2}. \]
Figure \ref{string-in-sphere-fig-2} shows the evolution of the curve through three cycles in the vibration; note that the shape of the wave varies dramatically and even becomes non-embedded at times.
\begin{figure}[h]
\includegraphics[width=1.8in]{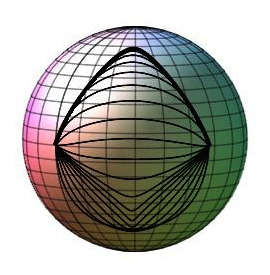}
\includegraphics[width=1.8in]{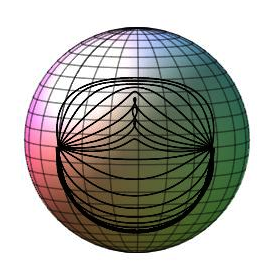}
\includegraphics[width=1.8in]{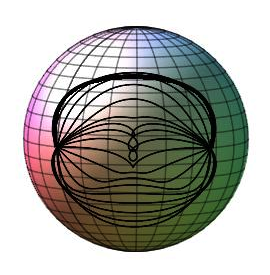}
\includegraphics[width=1.8in]{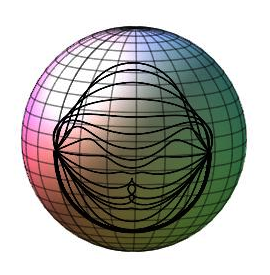}
\includegraphics[width=1.8in]{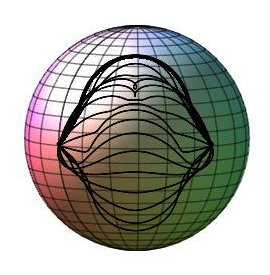}
\includegraphics[width=1.8in]{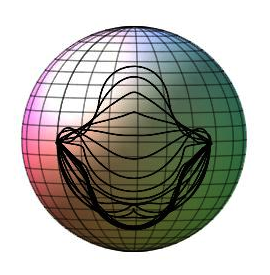}
\caption{Evolution of large amplitude string in $S^2$}
\label{string-in-sphere-fig-2}
\end{figure}

\end{example}

\begin{example}
Let $M = \mathbb{H}^2$ be the Poincar\'e upper half plane
\[ \mathbb{H}^2 = \{(x,y) \in \R^2 \ \vert \ y > 0 \}, \]
with its standard metric
\[  ds^2 = \frac{dx^2 + dy^2}{y^2}. \]
The wave map equation on $\mathbb{H}^2$ is equivalent to the system
\begin{align}
x_{tt} & = x_{mm} - \frac{2\,(x_m y_m - x_t y_t)}{y} , \label{H2-system} \\
y_{tt} & = y_{mm} + \frac{(x_m^2 - x_t^2) -  (y_m^2 - y_t^2)}{y} . \notag
\end{align}
As in the previous example, this is a nonlinear system, and its solutions are not periodic in general.

We will consider several different initial curves.  First consider the horizontal line segment $y=1$ for $0 \leq x \leq 2$, with parametrization
\[ x(m,0) = m, \qquad y(m,0) = 1, \qquad 0 \leq m \leq 2. \]
(Recall that this curve is not a geodesic in $\mathbb{H}^2$.)
Figure \ref{string-in-H2-fig-1} shows the evolution of the curve through three cycles in the vibration; the solution appears almost (but not quite) periodic.
\begin{figure}[h]
\includegraphics[width=1.8in]{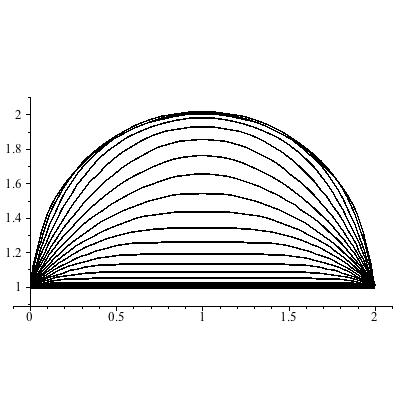}
\includegraphics[width=1.8in]{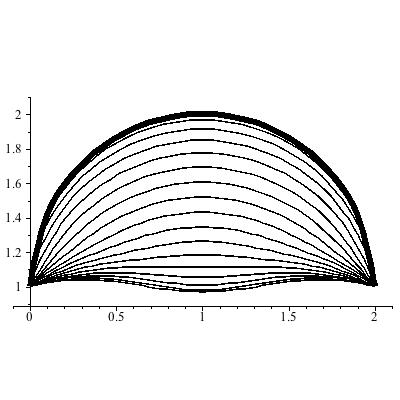}
\includegraphics[width=1.8in]{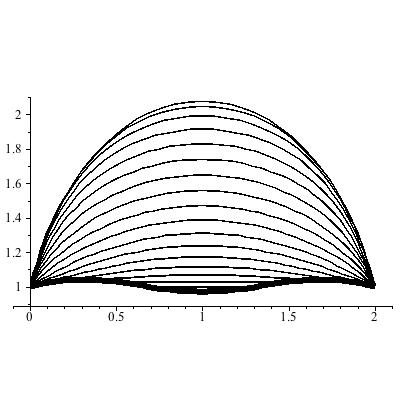}
\includegraphics[width=1.8in]{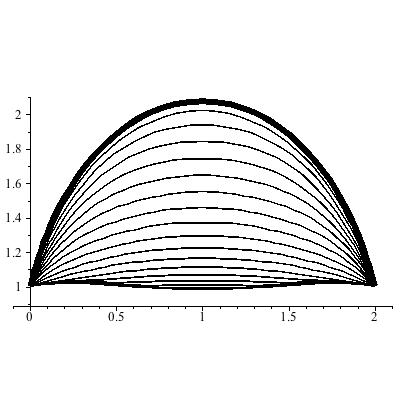}
\includegraphics[width=1.8in]{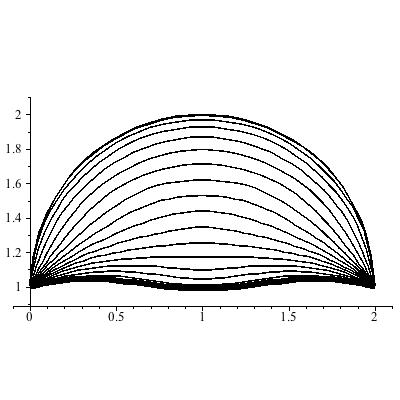}
\includegraphics[width=1.8in]{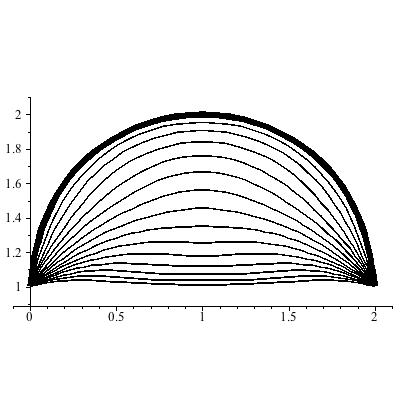}
\caption{Evolution of horizontal string in $\mathbb{H}^2$}
\label{string-in-H2-fig-1}
\end{figure}

Next, consider a ``sinusoidal" perturbation of a vertical geodesic segment.  In order that the initial curve appear approximately sinusoidal with respect to the Poincar\'e metric, we will take the curve
\[ x = \frac{1}{2} y \sin(y-1),\qquad  1 \leq y \leq \pi + 1, \]
as our initial curve, so that the horizontal line segment from the point $(x(y), y)$ to the point $(0,y)$ has length $\tfrac{1}{2} \sin(y-1)$.  (Note that this is still only approximately sinusoidal, since the horizontal line segment joining these two points is not a geodesic.)  Furthermore, we will parametrize this curve according to an arc length parametrization for the corresponding vertical line segment along the $y$-axis:
\[ x(m,0) = \frac{1}{2} e^m \sin(e^m-1),  \qquad  y(m,0) = e^m, \qquad 0 \leq m \leq \ln(\pi+1). \]
Figure \ref{string-in-H2-fig-2} shows the evolution of the curve through three cycles in the vibration; one intriguing feature of this evolution is that within each cycle, the string moves to the left much more rapidly than to the right.
\begin{figure}[h]
\includegraphics[width=1.8in]{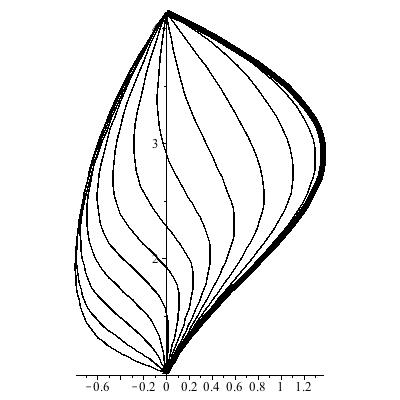}
\includegraphics[width=1.8in]{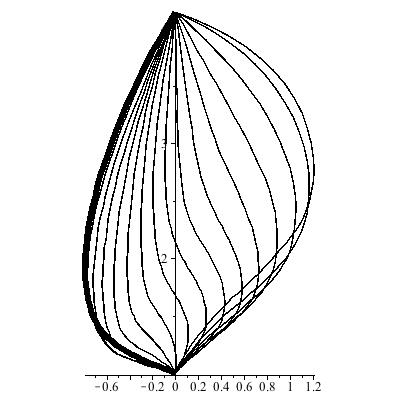}
\includegraphics[width=1.8in]{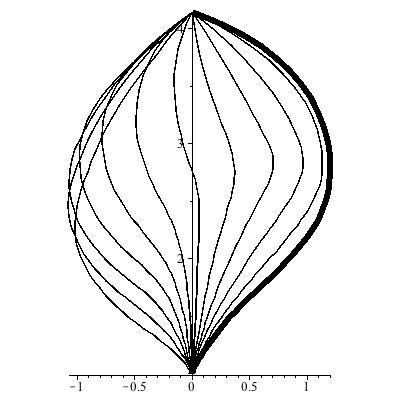}
\includegraphics[width=1.8in]{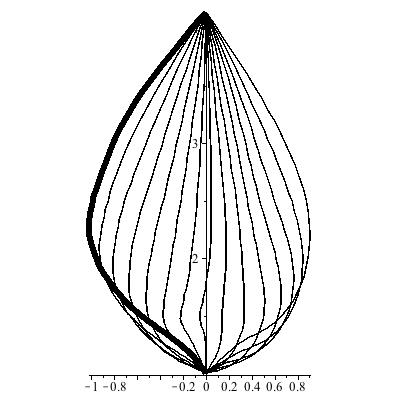}
\includegraphics[width=1.8in]{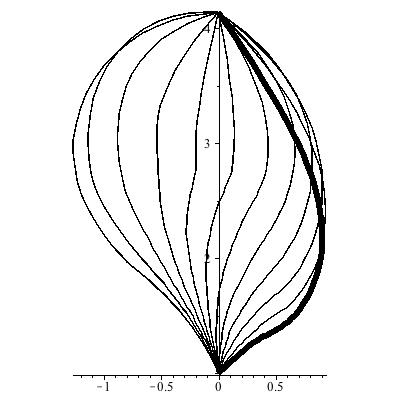}
\includegraphics[width=1.8in]{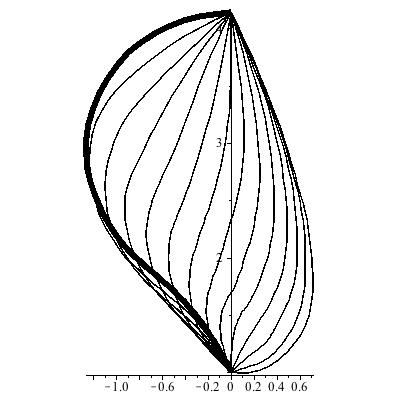}
\caption{Evolution of small amplitude string in $\mathbb{H}^2$}
\label{string-in-H2-fig-2}
\end{figure}

As in the previous example, the effect of the nonlinearity becomes more significant if we increase the amplitude of the initial curve: suppose that the initial curve is described by the equation
\[ x = y \sin(y-1), \qquad  1 \leq y \leq \pi + 1, \]
with parametrization
\[ x(m,0) = e^m \sin(e^m-1),  \qquad  y(m,0) = e^m, \qquad 0 \leq m \leq \ln(\pi+1). \]
Figure \ref{string-in-H2-fig-3} shows the evolution of the curve through three cycles in the vibration; now the shape of the wave varies dramatically and even becomes non-embedded at times.
\begin{figure}[h]
\includegraphics[width=1.8in]{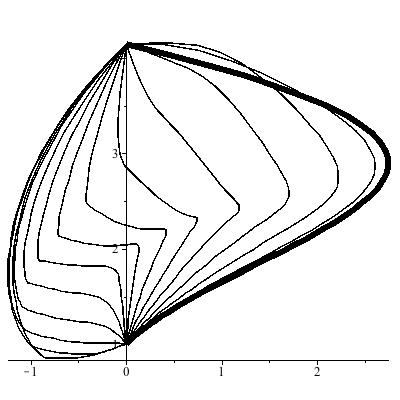}
\includegraphics[width=1.8in]{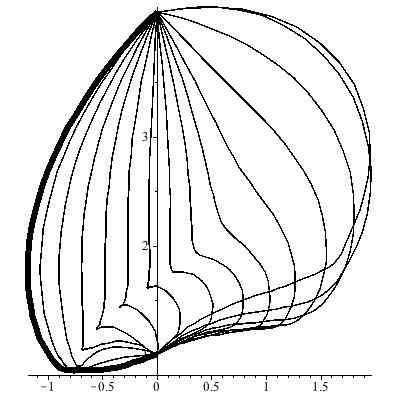}
\includegraphics[width=1.8in]{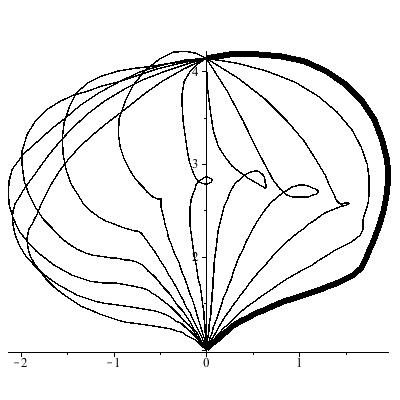}
\includegraphics[width=1.8in]{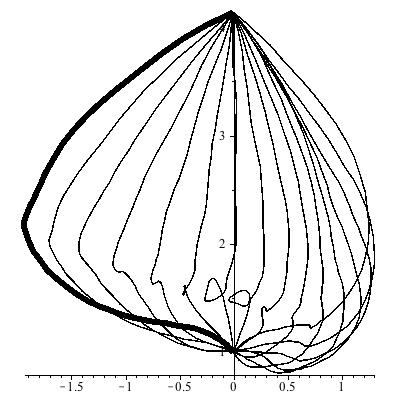}
\includegraphics[width=1.8in]{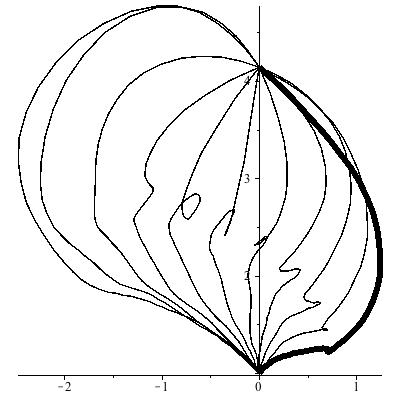}
\includegraphics[width=1.8in]{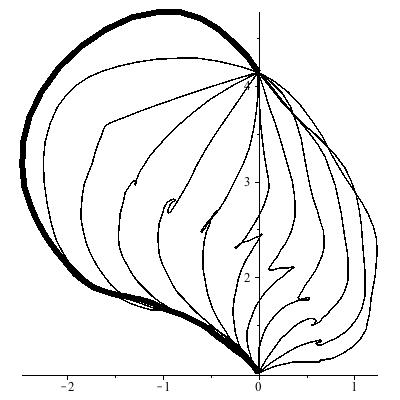}
\caption{Evolution of large amplitude string in $\mathbb{H}^2$}
\label{string-in-H2-fig-3}
\end{figure}

\end{example}

\section{Conclusion}

In this paper we have shown that the classical, frictionless, elastic string under tension undergoing transverse vibrations in the Euclidean plane is {\it exactly} modeled by the linear wave equation,
\begin{equation}\label{LinearWaveEqConclusion}
u_{tt}=c^2u_{xx},
\end{equation}
regardless of the size of the vibrations, and that this PDE is {\it not} a small amplitude approximation of a ``nonlinear string equation." More generally, if the arena of the vibration changes from the Euclidean plane to an arbitrary Riemannian surface, such the 2-sphere $S^2$ or the Poincare half-plane $\mathbb{H}^2$, then the motion of the string is governed by the partial differential equation for wave maps into Riemannian surfaces:
\begin{equation}\label{waveMapEqConclusion}
\nabla_{\bu_t} \bu_t = c^2 \nabla_{\bu_x} \bu_x.
\end{equation}
The wave map equation is intrinsically nonlinear whenever the Gauss curvature of the target space of the wave map is non-zero.  Equation (\ref{LinearWaveEqConclusion}) is essentially the wave map equation for wave maps into the Euclidean plane.  
Equation (\ref{waveMapEqConclusion}) is a source of fascination for both the mathematician and physicist. If forces cease to act on the string, then the left hand side of (\ref{waveMapEqConclusion}) will vanish, and the resulting differential equation will be that of a geodesic of the Riemannian surface.  Hence (\ref{waveMapEqConclusion}) generalizes the notion of geodesics and provides a very interesting class of hyperbolic partial differential equations for geometric analysis (\cite{Struwe}, \cite{TU04}). For the mathematician, there are still many basic questions that remain open. On the other hand, the same partial differential equation arises, for instance, in elementary particle physics (\cite{ketov}) and general relativity (\cite{Misner}, \cite{YCB1}, \cite{YCB2}).  

One natural question to consider is which Riemannian surfaces give rise to wave maps with particularly nice properties. In \cite{CV13} we have begun an exploration of wave maps for which the corresponding wave map system is {\em Darboux integrable}; we plan to delve more deeply into this topic in future work.  
Another area which we hope to explore in the near future is that of sub-Riemannian geometry; the context in which we have discussed wave maps here as parametrized curves evolving in time still makes sense in sub-Riemannian geometry, and questions such as global existence of solutions and integrability are of interest.

\nocite{Vassiliou12}

\bibliographystyle{amsplain}
\bibliography{StringsAttached-bib}

\end{document}